\documentclass[pdflatex,12pt]{elsarticle}

\usepackage{amssymb}
\usepackage{amsmath}
\usepackage{hyperref}
\hypersetup{hypertex=true,
            colorlinks=true,
            linkcolor=blue,
            anchorcolor=blue,
            citecolor=blue}
\usepackage{amsthm}
\usepackage{graphicx}%
\usepackage{multirow}%
\usepackage{amsmath,amssymb,amsfonts}%
\usepackage{amsthm}%
\usepackage{mathrsfs}%
\usepackage[title]{appendix}%
\usepackage{xcolor}%
\usepackage{textcomp}%
\usepackage{manyfoot}%
\usepackage{booktabs}%
\usepackage{algorithm}%
\usepackage{algorithmicx}%
\usepackage{algpseudocode}%
\usepackage{tcolorbox}  
\tcbuselibrary{most}     
\newtcolorbox{examplebox}{ 
    colback=white,          
    colframe=black!80!white,
    boxrule=0.8pt,          
    arc=0pt,                
    left=4mm, right=4mm,    
    top=2mm, bottom=2mm,    
    before upper={\parindent0pt\noindent}, 
    breakable,              
    enhanced                
}

\usepackage{listings}%
\theoremstyle{thmstyleone}%
\newtheorem{theorem}{Theorem}
\theoremstyle{thmstyletwo}%
\newtheorem*{remark*}{Remark}
\theoremstyle{thmstylethree}%
\newtheorem{lemma}{Lemma}

\newtheorem{cor}{Corollary}
\newtheorem{conje}{Conjecture}
\newtheorem*{example*}{Example}
\journal{European Journal of Combinatorics}

\begin{document}

\begin{frontmatter}
\title{On the biases and asymptotics of partitions with finite choices of parts}

\author{Jiyou Li}
 \ead{lijiyou@sjtu.edu.cn}
\affiliation{organization={School of Mathematical Sciences, Shanghai Jiao Tong University},
            city={Shanghai},
            postcode={200240},
            country={China}
             }

\author{Sicheng Zhao} 
\ead{zsc_abandoned@alumni.sjtu.edu.cn}
\affiliation{organization={School of Mathematical Sciences, Shanghai Jiao Tong University},
            city={Shanghai},
            postcode={200240},
            country={China}
            }


\abstract{Biases in integer partitions have been studied recently. For three disjoint subsets $R,S,I$ of positive integers, let $p_{RSI}(n)$ be the number of partitions of $n$ with parts from $R\cup S\cup I$ and $p_{R>S,I}(n)$ be the number of such partitions with  a greater number of parts in $R$ than that in $S$. In this paper,  in the case that $R,S,I$ are finite, we obtain an explicit formula of the asymptotic ratio of $p_{R>S,I}(n)$ to $p_{RSI}(n)$. The key technique for computing this ratio is to estimate a partition number at the volume of a certain polytope. A conjecture is proposed in the case that $R,S$ are certain infinite  arithmetic
progressions.}

}
\end{frontmatter}

\maketitle

\section{Introduction}\label{sec1}
Partition theory is one of the most important subjects in combinatorics. For a positive integer $n$, a partition of $n$ is a non-increasing positive integer sequence $\lambda=(\lambda_1,\lambda_2,\ldots,\lambda_l)$ with $\sum\limits_{i=1}^l \lambda_l=n$. Each $\lambda_i$ is called a part and $l$ is called the length of $\lambda$. 
Let $p(n)$ denote the number of partitions of $n$. One of the most remarkable results about $p(n)$ is its asymptotic estimate given by Hardy and Ramanujan \cite{Ram}:
\begin{equation*}
    p(n)\sim \frac{e^{\pi\sqrt{\frac{2n}{3}}}}{4\sqrt{3}n}.
\end{equation*}
This is actually the origin of the circle method which plays an important role in combinatorics and number theory.

Partitions with some constraints have been studied by many mathematicians since Euler. Many related questions also arise naturally from algebra, such as the representation theory of symmetric groups or the theory of modular forms. Extensive research has been carried out on the study of partitions with restricted parts in many different kind of settings.

A typical example of restricted partitions is the $t$-regular partition, that is, the partition whose parts are not divisible by $t$. The study of the $t$-regular partitions can be traced back to Euler \cite{Euler} and Glaisher \cite{Gla} several centuries ago. The arithmetic properties of the number of $t$-regular partitions of $n$, such as congruence or divisibility, are also of great interest to researchers. For example, see \cite{MR2501237}, \cite{MR4632990} and \cite{MR4707380}.

Recently ‘‘biases in partition’’ arise as an interesting subject in the study of partitions. By biases in partition, we mean the tendency of some parts to appear more frequently than other parts in the partitions (may be restricted). Quantitative and qualitative researches have been carried out on such phenomena. 

For instance,  the biases in the appearance of the parts from two residue classes in the partitions have been studied recently. Kim et al. 
\cite{Kim1} investigate the number of partitions of $n$ with more (or less) odd parts than even parts, which will be denoted as $p_o(n)$ (or $p_e(n)$). They prove that $p_o(n)>p_e(n)$ for $n$ large enough and name this phenomenon as parity bias. Precisely, the authors show that
\begin{equation*}
    \begin{aligned}
        p_o(n)&\sim \frac{1}{\sqrt{2}}p(n),\\
        p_e(n)&\sim (1-\frac{1}{\sqrt{2}})p(n).
    \end{aligned}
\end{equation*}

In a following study \cite{Kim2} of Kim and Kim, the biases between certain residue classes are studied. Let $p_{a,b,m}(n)$ be  the number of partitions of $n$ with more parts congruent to $a$ modulo $m$ than parts congruent to $b$ modulo $m$ for $m\geq 2$. They prove that 
\begin{equation}\label{1mm}
    \begin{aligned}
        p_{1,m,m}(n)&\sim \frac{1}{2^{\frac{1}{m}}}p(n),\\
        p_{m,1,m}(n)&\sim (1-\frac{1}{2^{\frac{1}{m}}})p(n).
    \end{aligned}
\end{equation}
This result indicates that the partitions of $n$ tend to have more  parts congruent to $1$ modulo $m$ than parts congruent to $m$ modulo $m$ for $n$ large enough. In general, Chern \cite{Che} shows that for $1\leq a<b\leq m$,
\begin{equation*}
    p_{a,b,m}(n)\geq p_{b,a,m}(n).
\end{equation*}
Interested readers are also referred to \cite{D}, \cite{B}, \cite{nonuni}, \cite{Kim3}, \cite{kathrin}, \cite{Kim4} for more such results.

It is natural to note that among all partitions of a positive integer $n$,  the  smaller parts appear more frequently than the larger parts. However, it seems nontrivial to quantitatively characterize such phenomenon, which is the motivation of this paper.

The main purpose of this paper is to explicitly measure the biases of the appearance of parts from two certain finite sets when the parts of partitions have only finite choices. Let $R,S,I$  be three disjoint subsets of $\mathbb{Z}^+$. Let $P_{RSI}(n)$ denote the set of partitions of $n$ with parts in $R\cup S\cup I$ and $P_{R>S,I}(n)$ denote the set of partitions in $P_{RSI}(n)$ with a greater number of  parts in $R$ than that in $S$. Throughout this paper,  the lowercase letters $p_{RSI}(n),p_{R>S,I}(n)$ will represent the size of a set of partitions. The main result of this paper is an interesting formula given in the following theorem. 
~\\
\begin{theorem}\label{thm1}
Let $R=\{r_1,\ldots,r_l\},S=\{s_1,s_2\ldots,s_m\}$ with $R,S$ be disjoint positive integer sets and $I$ be a finite subset of $\mathbb{Z}^+\setminus(R\cup S)$. Suppose the greatest common divisor of $\{r_j\}_{j=1}^l,\{s_i\}_{i=1}^m$ is $1$. Then
\begin{equation}\label{eq3}
    \lim\limits_{n\rightarrow \infty}\frac{p_{R>S,I}(n)}{p_{RSI}(n)}=\prod\limits_{j=1}^l r_j\prod\limits_{i=1}^m s_i\cdot\sum\limits_{i=1}^l\frac{(-1)^{i-1}}{r_i\prod\limits_{j=i+1}^l(r_j-r_i)\prod\limits_{t=1}^{i-1}(r_i-r_t)\prod\limits_{k=1}^m(s_k+r_i)}.
\end{equation}
\end{theorem}

The key technique for the proof of Theorem \ref{thm1} is to estimate the partition numbers at the volume of certain polytope. Firstly, we reduce the condition to the case that $I$ is empty. Then we establish an one-to-one correspondence between the partitions in  $P_{RSI}(n)$ (or $P_{R>S,I}(n)$) and the integer points in a certain finite dimensional polytope. Finally, the asymptotic estimates of $p_{RSI}(n)$ and $p_{R>S,I}(n)$ are given by some technical but elementary calculations of the integral on the polytopes.

Several interesting corollaries are listed here as examples of the applications of Theorem \ref{thm1}. The proofs are direct and hence omitted. 
~\\
\begin{cor}

Let $R=\{r\}$ and $S=\{s\}$ with $(r,s)=1$, $I$ be an  finite subset of $\mathbb{Z}^+\setminus(R\cup S)$. Then
\[\lim\limits_{n\rightarrow \infty}\frac{p_{R>S,I}(n)}{p_{RSI}(n)}=\frac{s}{r+s}.\]
\end{cor}

\begin{cor}
Let $R=\{1\}$ and $S=\{s_1,s_2,\ldots,s_m\}$ with $s_i>1$, $I$ be an  finite subset of $\mathbb{Z}^+\setminus(R\cup S)$. Then   
\[\lim\limits_{n\rightarrow \infty}\frac{p_{R>S,I}(n)}{p_{RSI}(n)}=\prod\limits_{i=1}^m \frac{s_i}{s_i+1}.\]
In particular, if $S=\{2, 3, \ldots, k\}$ with $k\geq 2$, one has
\[\lim\limits_{n\rightarrow \infty}\frac{p_{R>S,I}(n)}{p_{RSI}(n)}=\frac 2 {k+1}.\]
\end{cor}

\begin{cor}
Let  $R=\{1,2\}$ and $S=\{3, 4\ldots,k\}$ with $k\geq 3$, $I$ be a finite subset of $\mathbb{Z}^+\setminus(R\cup S)$. Then
\[\lim\limits_{n\rightarrow \infty}\frac{p_{R>S,I}(n)}{p_{RSI}(n)}=\frac{6k}{(k+1)(k+2)}.\]
\end{cor}
~\\

Theorem \ref{thm1} gives the asymptotic proportion of the  partitions with more parts from $R$ than parts from $S$ in $P_{RSI}(n)$. Note that  when the RHS of (\ref{eq3}) is greater than $\dfrac{1}{2}$, partitions in $P_{RSI}(n)$ tend to have more parts from $R$ than parts from $S$.

It is very natural to consider the case when $R, S$ or $I$ are infinite. Building upon the results of \cite{Kim2} , we propose a conjecture on the infinite cases, motivated by the applications of our formula (\ref{eq3}) in certain finite cases. 

Assume  $r\not\equiv s\ (mod\ m)$ and $(r,s,m)=1$. Let $R_N=\{r,r+m,\ldots,r+m(N-1)\}$ and $S_N=\{s,s+m,\ldots,s+m(N-1)\}$ be two arithmetic progressions of the same length but with distinct starting number, and let $I_N=\left[\max(r,s)+m(N-1)\right]\setminus(R_N\cup S_N)$. Let $R=\{r,r+m,\ldots,r+im,\ldots\},S=\{s,s+m,\ldots,s+im,\ldots\}, I=\mathbb{Z}^+\setminus(R\cup S)$. Let $C_{n,N}=\dfrac{p_{R_N>S_N,I_N}(n)}{p_{R_NS_NI_N}(n)}$.  Note that for a fixed $n$,  $C_{n,N}=\dfrac{p_{R>S,I}(n)}{p_{RSI}(n)}$ when $N$ is large enough. That is,
\begin{equation*}
    \lim\limits_{n\rightarrow \infty}\frac{p_{R>S,I}(n)}{p_{RSI}(n)}=\lim\limits_{n\rightarrow \infty}\lim\limits_{N\rightarrow \infty}C_{n,N}.
\end{equation*}

We conjecture that

\begin{equation*}
        \lim\limits_{n\rightarrow \infty}\lim\limits_{N\rightarrow \infty}C_{n,N}=\lim\limits_{N\rightarrow \infty}\lim\limits_{n\rightarrow \infty}C_{n,N}.
\end{equation*}

Therefore, if the conjecture is true, we can evaluate $\lim\limits_{N\rightarrow \infty}\lim\limits_{n\rightarrow \infty}C_{n,N}$ instead of $ \lim\limits_{n\rightarrow \infty}\dfrac{p_{R>S,I}(n)}{p_{RSI}(n)}$. Then Theorem \ref{thm1} can directly apply to the computation of $\lim\limits_{N\rightarrow \infty}\lim\limits_{n\rightarrow \infty}C_{n,N}$. Actually, the conjecture is motivated by the numerical agreement of the computational result from Theorem \ref{thm1} and that derived in \cite{Kim2} when $(r,s,m)=(1,m,m)$. We direct readers to more detailed discussion in Section \ref{sec3}.

The rest of this paper is organized as follows. The proof of Theorem \ref{thm1} is given in Section \ref{sec2} by applying area integral to estimate the partition numbers.  Further discussion of the applications of Theorem \ref{thm1} and the above conjecture is presented in Section \ref{sec3}.

\section{Proof of Theorem \ref{thm1}}\label{sec2}

In the following lemma, we first show  that $I$ is independent to the asymptotic value of $\dfrac{P_{R>S,I}(n)}{P_{RSI}(n)}$ if it is finite.
~\\
\begin{lemma}\label{lem1}
Let $R,S,I$ and other notations be defined above. Suppose $i\in \mathbb{Z}^+\setminus{(R\cup S\cup I)}$ and 

    \begin{equation*}
    \lim\limits_{n\rightarrow \infty}\frac{p_{R>S,I}(n)}{p_{RSI}(n)}=C.
\end{equation*}
Then for $I'=I\cup\{i\}$,
\begin{equation*}
    \lim\limits_{n\rightarrow \infty}\frac{p_{R>S,I'}(n)}{p_{RSI'}(n)}=C.
\end{equation*}
\end{lemma}
~\\
\begin{proof}
    Let 
    \begin{equation*}
    \begin{aligned}
        f(q)&=\sum\limits_{n=1}^\infty p_{R>S,I}(n)q^i,\\
        g(q)&=\sum\limits_{n=1}^\infty p_{RSI}(n)q^i,
    \end{aligned}
    \end{equation*}
    be the generating functions of ${p_{R>S,I}(n)},\ {p_{RSI}(n)}$. Then clearly the generating functions of ${p_{R>S,I'}(n)},\ {p_{RSI'}(n)}$ are $\dfrac{f(q)}{1-q^i},\ \dfrac{g(q)}{1-q^i}$. Therefore
    \begin{equation*}
        \frac{p_{R>S,I'}(n)}{p_{RSI'}(n)}=\frac{p_{R>S,I}(n)+p_{R>S,I}(n-i)+p_{R>S,I}(n-2i)+\ldots}{p_{RSI}(n)+p_{RSI}(n-i)+p_{RSI}(n-2i)+\ldots}.
    \end{equation*}
Separate the sequence $\left\{\dfrac{p_{R>S,I'}(n)}{p_{RSI'}(n)}\right\}_{n=1}^{\infty}$ into $i$ disjoint subsequences by the residue of $n$ module $i$ and apply the Stolz-Cesaro theorem (Theorem 1.22, \cite{stolzthm}). The lemma then follows directly.

\end{proof}

By lemma \ref{lem1}, we may assume $I$ to be empty in the proof of Theorem \ref{thm1}. For simplicity, $P_{RSI}(n),P_{R>S,I}(n)$ will be written as $P_{RS}(n),P_{R>S}(n)$ in the rest of this section. We use multisets to denote partitions to avoid extra notation on the order of the parts. The superscript will represent the multiplicity of the number as a part in a multiset (namely, a partition). For example, $\lambda=(4,4,2,1,1,1)$, a partition of $13$, will be written as $\{4^2,2^1,1^3\}$.

Now let us first consider $P_{RS}(n)$. Since the greatest common divisor of $\{r_j\}_{j=1}^l,\{s_i\}_{i=1}^m$ is $1$, let
\begin{equation*}
    \sum\limits_{j=1}^l a_jr_j+\sum\limits_{i=1}^m b_is_i=1,
\end{equation*}
where $\{a_j\}_{j=1}^l,\{b_i\}_{i=1}^m$ are integers. Then
\begin{equation*}
    \sum\limits_{j=1}^l na_j\cdot r_j+\sum\limits_{i=1}^m nb_i\cdot s_i=n.
\end{equation*}
$\forall\lambda=\{r_1^{c_1},\ldots,r_l^{c_l},s_1^{f_1},\ldots,s_m^{f_m}\}\in P_{RS}(n)$, we have
\begin{equation*}
    \sum\limits_{j=1}^l (na_j-c_j)\cdot r_j+\sum\limits_{i=1}^m (nb_i-f_i)\cdot s_i=0.
\end{equation*}

Let $na_j-c_j=x_j$ for $1\leq j\leq l$ and $nb_i-f_i=y_i$ for $1\leq i\leq m$.  Then $\{x_j\}_{j=1}^l,\{y_i\}_{i=1}^m$ are integers. We have
\begin{equation}\label{lat1}
     \sum\limits_{j=1}^l x_j\cdot r_j+\sum\limits_{i=1}^m y_i\cdot s_i=0.
\end{equation}

(\ref{lat1}) is actually an affine hyperplane in $\mathbb{R}^{l+m}$ and its integer points form a lattice $\Lambda_1$ of dimension $l+m-1$. Then $p_{RS}(n)$ equals the number of lattice points of $\Lambda_1$ such that the corresponding $\{c_j\}_{j=1}^l,\{f_i\}_{i=1}^m$ are non-negative. 
\vspace{0.7cm}
  \begin{example*}
   Let $R=\{2,3,6\},S=\{10,15\}$. Let $(a_1,a_2,a_3,b_1,b_2)=(-1,1,0,0,0)$. For a partition $\{2^{c_1},3^{c_2},6^{c_3},10^{f_1},15^{f_2}\}\in P_{RS}(n)$, we have
   \begin{equation*}
       2c_1+3c_2+6c_3+10f_1+15f_2=n,
   \end{equation*}
   
Therefore, for $x_1=-n-c_1,x_2=n-c_2,x_3=-c_3,y_1=-f_1,y_2=-f_2$, $\{x_j\}_{j=1}^l,\{y_i\}_{i=1}^m$ are integers satisfying the following equation:
\begin{equation*}
    2x_1+3x_2+6x_3+10y_1+15y_2=0.
\end{equation*}
The integer points of the above equation form a lattice of dimension $4$ in $\mathbb{R}^{5}$. Any lattice point satisfying $x_1\leq-n,\ x_2\leq n,\ x_3\leq 0,\ y_1\leq 0,\ y_2\leq 0$ is in one-to-one correspondence with a partition in $P_{RS}(n)$.
\end{example*}  
\vspace{0.7cm}

The above example will serve as a running example throughout this section.

To estimate the number of such points, let us first construct a basis of this lattice.

Since $r_1$ is coprime to $(r_2,\ldots,r_l,s_1,\ldots,s_m)$, the minimal possible positive value of $x_1$ is $d_1=(r_2,\ldots,r_l,s_1,\ldots,s_m)$. Select an arbitrary vector $\mathbf{v}_1$ of $\Lambda_1$ with $x_1=d_1$. Then the quotient lattice $\Lambda_2=\Lambda_1/\left\langle \mathbf{v}_1\right\rangle$ is a lattice of dimension $l+m-2$ determined by following equations:
\begin{equation*}
    x_1=0,\ \sum\limits_{j=2}^l x_j\cdot r_j+\sum\limits_{i=1}^m y_i\cdot s_i=0.
\end{equation*}

In $\Lambda_2$, the minimal possible positive value of $x_2$ is $d_2=\dfrac{(r_3,\ldots,r_l,s_1,\ldots,s_m)}{(r_2,\ldots,r_l,s_1,\ldots,s_m)}$. Select an arbitrary vector $\mathbf{v}_2$ of $\Lambda_2$ with $x_2=d_2$. Then the quotient lattice $\Lambda_3=\Lambda_2/\left\langle \mathbf{v}_2\right\rangle$ is a lattice of dimension $l+m-3$ determined by following equations:
\begin{equation*}
    x_1=x_2=0,\ \sum\limits_{j=3}^l x_j\cdot r_j+\sum\limits_{i=1}^m y_i\cdot s_i=0.
\end{equation*}

Repeat this operation until we have already selected $l+m-1$ vectors $\mathbf{v}_1,\mathbf{v}_2,\ldots,\mathbf{v}_{l+m-1}$.  If we rename $r_1,\ldots,r_l,s_1,\ldots,s_m$ as $e_1,\ldots,e_{l+m}$, then the first non-zero component of $\mathbf{v}_i$ is on the $i$-th position and its value is $\dfrac{(e_{i+1},e_{i+2},\ldots,e_{l+m})}{(e_i,e_{i+1},\ldots,e_{l+m})}$. One can easily check that $\{\mathbf{v}_1,\mathbf{v}_2,\ldots,\mathbf{v}_{l+m-1}\}$ form a basis of $\Lambda_1$.
~\\
\vspace{0.7cm}
\begin{example*}
    Let $R=\{2,3,6\},S=\{10,15\}$. We could select $\mathbf{v}_1$ to be $(1,-4,0,1,0)$, and then $\mathbf{v}_2=(0,1,-3,0,1),\ \mathbf{v}_3=(0,0,5,-3,0),\ \mathbf{v}_4=(0,0,0,3,-2)$. So the row vectors of the following matrix form a basis of $\Lambda_1$:
\begin{equation*}
\begin{pmatrix}
    1&-4&0&1&0  \\
     0&1&-3&0&1 \\
     0&0&5&-3&0\\
     0&0&0&3&-2\\
\end{pmatrix}
\end{equation*}

\end{example*}
\vspace{0.7cm}

Now let us return to the estimate of the number of the lattice points of $\Lambda_1$ such that the corresponding $\{c_j\}_{j=1}^l,\{f_i\}_{i=1}^m$ are non-negative. Let $\mathbf{w}_n=(na_1,\ldots,na_l,nb_1,\ldots,nb_m)$ and $(x_1,\ldots,x_l,y_1,\ldots,y_m)=-\sum\limits_{i=1}^{l+m-1}k_i\mathbf{v}_i$, then
\begin{equation}\label{sys1}
\left.
\begin{aligned}
   (c_1,\ldots,c_l,f_1,\ldots,f_m)&=\mathbf{w}_n-(x_1,\ldots,x_l,y_1,\ldots,y_m)\\
   &=\mathbf{w}_n+\sum\limits_{i=1}^{l+m-1}k_i\mathbf{v}_i \\
   &\geq \mathbf{0}.
\end{aligned}
\right.
\end{equation}

Here "$\geq$" is component-wise. Then $p_{RS}(n)$ is equal to the number of integer arrays $\{k_i\}_{i=1}^{l+m-1}$ satisfying  system (\ref{sys1}), namely the number of integer points in the area determined by (\ref{sys1}) in variables $\{k_i\}_{i=1}^{l+m-1}$. 
\vspace{0.7cm}
    \begin{example*}
    Let $R=\{2,3,6\},S=\{10,15\}$ and let $a_i,b_i,\mathbf{v}_i$'s be defined as previous examples. Then the following system derived from (\ref{sys1})
\begin{equation*}
\left\{
\begin{aligned}
k_1&&&&&\geq n,\\
-4k_1&+\ k_2&&&&\geq -n,\\
&-3k_2&+5k_3&&&\geq 0,\\
k_1&&-3k_3&+3k_4&&\geq 0,\\
&\ \ \ \ \ k_2&&-2k_4&&\geq 0,\\
\end{aligned}
\right.
\end{equation*}
determines an area of $\mathbb{R}^{4}$. Any integer array $\{k_i\}_{i=1}^{4}$ satisfying  this system is in bijective correspondence with a partition in $P_{RS}(n)$. 

For instance, set $n=10$ and let $\{2^5\}\in P_{RS}(10)$. By calculation, the corresponding $(c_1,c_2,c_3,f_1,f_2)=(5,0,0,0,0)$, $(x_1,x_2,x_3,y_1,y_2)=(-15,10,0,0,0)$. Then $(k_1,k_2,k_3,k_4)=(15,50,30,25)$, which satisfies the system above.

\end{example*}
\vspace{0.7cm}

Denote this area determined by (\ref{sys1}) in variables $\{k_i\}_{i=1}^{l+m-1}$ as $D_n$. Note that $D_n$ is actually  the intersection of $l+m$ half-spaces in $\mathbb{R}^{l+m-1}$, so it is a polytope of dimension $l+m-1$. Furthermore, the $l+m$ inequations of (\ref{sys1}) are linear in variables $\{k_i\}_{i=1}^{l+m-1}$ and the constant terms are proportional to $n$, namely $D_n$ is the $n$-dilate of $D_1$. It is time to introduce the following lemma from  Ehrhart theory.
~\\
\begin{lemma}[Lemma 3.19, \cite{MaSi}]\label{lem2}
    Suppose $P\subseteq \mathbb{R}^d$ is $d$-dimensional. Then
\begin{equation*}
    vol(P)=\lim\limits_{t\rightarrow \infty}\frac{1}{t^d}\cdot\# (tP\cap \mathbb{Z}^d). 
\end{equation*}
    
\end{lemma}
Therefore an asymptotic estimate of the number of integer points in  $D_n$, namely $\# (D_n\cap \mathbb{Z}^d)$, is $vol(D_1)\cdot n^{l+m-1}+o(n^{l+m-1})$, where $vol(D_1)=\int_{R^{l+m-1}}\chi_{D_1}\ dV$. Let $u_1 =c_1,\ldots,u_l=c_l,u_{l+1}=f_1,\ldots,u_{l+m-1}=f_{m-1}$ and set $n=1$. Then (\ref{sys1}) is transformed into
\begin{equation}\label{sys2}
\left\{
\begin{aligned}
&u_i\geq 0,\ 1\leq i\leq l+m-1,\\
&r_1u_1+\ldots+r_lu_l+s_1u_{l+1}+\ldots+s_{m-1}u_{l+m-1}\leq 1.\\
\end{aligned}
\right.
\end{equation}

Here the last inequality comes from $\sum\limits_{j=1}^l r_j c_j+\sum\limits_{i=1}^m s_i f_i=n$ and $f_m\geq 0$. Due to the discussion above, the elements in the main diagonal of the Jacobi matrix of this variable substitution from $\{k_i\}_{i=1}^{l+m-1}$ to $\{u_i\}_{i=1}^{l+m-1}$
 are $\dfrac{(e_{i+1},e_{i+2},\ldots,e_{l+m})}{(e_i,e_{i+1},\ldots,e_{l+m})}(1\leq i\leq l+m-1)$ and its determinant is $e_{l+m}=s_m$. So the volume of $D_1$ is $\dfrac{1}{s_m}$ as the volume of (\ref{sys2}). Let $y_i=r_iu_i$ for $1\leq i\leq l$ and $y_i=s_{i-l}u_i$ for $l+1\leq i\leq l+m-1$. Then 

\begin{equation*}
    \left.
    \begin{aligned}
        \int_{\mathbb{R}^{l+m-1}}\chi_{D_1}\ dV&=\frac{1}{s_m}\frac{1}{\prod\limits_{i=1}^{m-1} s_i\prod\limits_{j=1}^l r_j}\int_0^{1}dy_1\int_0^{1-y_1}dy_2\ldots\int_0^{1-y_1-\ldots-y_{l+m-2}}  dy_{l+m-1}\\
        &=\frac{1}{(l+m-1)!}\frac{1}{\prod\limits_{i=1}^m s_i \prod\limits_{j=1}^l r_j}.
    \end{aligned}
    \right.
\end{equation*}
Here the second equality comes from an easy computation of multiple integral. Therefore
\begin{equation}\label{eq1}
    p_{RS}(n)=\frac{A_{l+m-1}}{\prod\limits_{i=1}^m s_i \prod\limits_{j=1}^l r_j}+o(n^{l+m-1}),
\end{equation}
where $A_{l+m-1}$ denotes $\dfrac{n^{l+m-1}}{(l+m-1)!}$.

To estimate $p_{R>S}(n)$, add a new condition $\sum\limits_{j=1}^l c_j>\sum\limits_{i=1}^m f_i$ to system (\ref{sys1}). Note that the argument above still applies here. Applying the same variable substitution and setting $n=1$, we have
\begin{equation}\label{sys3}
\left\{
\begin{aligned}
&u_i\geq 0, 1\leq i\leq l+m-1,\\
&r_1u_1+\ldots+r_lu_l+s_1u_{l+1}+\ldots+s_{m-1}u_{l+m-1}\leq 1,\\
&\sum\limits_{j=1}^l(s_m+r_j)u_j+\sum\limits_{i=1}^{m-1}(s_i-s_m)u_{l+i}>1.\\
\end{aligned}
\right.
\end{equation}

Here the last inequality comes from $\sum\limits_{j=1}^l c_j>\sum\limits_{i=1}^m f_i$ and $\sum\limits_{j=1}^l r_j c_j+\sum\limits_{i=1}^m s_i f_i=n$. Denote the volume of  system (\ref{sys3}) as $V$. Then $p_{R>S}(n)=\dfrac{1}{s_m}\cdot V\cdot n^{l+m-1}+o(n^{l+m-1})$. For the ease of computation below, we will calculate the volume of the complement of (\ref{sys3}) in (\ref{sys2}) instead, namely the following system:
\begin{equation}\label{sys4}
    \left\{
\begin{aligned}
&u_i\geq 0,\ 1\leq i\leq l+m-1,\\
&r_1u_1+\ldots+r_lu_l+s_1u_{l+1}+\ldots+s_{m-1}u_{l+m-1}\leq 1,\\
&\sum\limits_{j=1}^l(s_m+r_j)u_j+\sum\limits_{i=1}^{m-1}(s_i-s_m)u_{l+i}\leq 1.\\
\end{aligned}
\right.
\end{equation}

Then the volume of (\ref{sys4}) equals $\dfrac{1}{(l+m-1)!\prod\limits_{j=1}^l r_j\prod\limits_{i=1}^{m-1} s_i}-V$. Separate (\ref{sys4}) into following two systems:
\begin{equation}\label{sys5}
\left\{
\begin{aligned}
&u_i\geq 0,\ 1\leq i\leq l+m-1,\\
&r_1u_1+\ldots+r_lu_l+s_1u_{l+1}+\ldots+s_{m-1}u_{l+m-1}\leq 1,\\
&\sum\limits_{j=1}^l u_j<\sum\limits_{i=1}^{m-1}u_{l+i}.\\
\end{aligned}
\right.
\end{equation}
and
\begin{equation}\label{sys6}
\left\{
\begin{aligned}
&u_i\geq 0,\ 1\leq i\leq l+m-1,\\
&\sum\limits_{j=1}^l(s_m+r_j)u_j+\sum\limits_{i=1}^{m-1}(s_i-s_m)u_{l+i}\leq 1,\\
&\sum\limits_{j=1}^l u_j\geq\sum\limits_{i=1}^{m-1}u_{l+i}.\\
\end{aligned}
\right.
\end{equation}

Note that (\ref{sys5}),(\ref{sys6}) are of the same type: if we define (\ref{sys5}) as an $(l,m-1)$-form, then (\ref{sys6}) is actually an $(m-1,l)$-form. Denote the volume of (\ref{sys5}) as $V_{(r_1,\ldots,r_l),(s_1,\ldots,s_{m-1})}$, then the volume of (\ref{sys6}) can be written as $V_{(s_1-s_m,\ldots,s_{m-1}-s_m),(r_1+s_m,\ldots,r_l+s_m)}$. These two volumes will be written as $V_1,V_2$ for convenience sometime below. Now let us compute the $(l,m-1)$-form (\ref{sys5}) .

Let $v=\sum\limits_{i=1}^{m-1}u_{l+i}-\sum\limits_{j=1}^l u_j$ and eliminate $u_{l+1}$, (\ref{sys5}) is transformed into:
\begin{equation}\label{sys7}
\left\{
\begin{aligned}
&u_i\geq 0,\ 1\leq i\leq l+m-1,\ i\neq l+1,\\
&v\geq 0,\\
&(s_2-s_1)u_{l+2}+\ldots+(s_{m-1}-s_1)u_{l+m-1}+s_1v+(r_1+s_1)u_1+\ldots+(r_l+s_1)u_l\leq 1,\\
&\sum\limits_{i=2}^{m-1}u_{l+i}\leq v+\sum\limits_{j=1}^l u_j.\\
\end{aligned}
\right.
\end{equation}

Since the Jacobi matrix of this transformation is upper-triangular and the diagonal elements are $1$, the volume of ($\ref{sys5}$)  is equal to the volume of ($\ref{sys7}$), namely
\begin{equation}\label{relation2}
    V_{(r_1,\ldots,r_l),(s_1,\ldots,s_{m-1})}=V_{(s_2-s_1,\ldots,s_{m-1}-s_1),(s_1,r_1+s_1,\ldots,r_l+s_1)}.
\end{equation}
So the computation of an $(l,m-1)$-form is equivalent to that of an $(m-2,l+1)$-form.
\vspace{0.7cm}
\begin{example*}
 Let $R=\{2,3,6\},S=\{10,15\}$, then system (\ref{sys4}) is transformed to
 \begin{equation*}
\left\{
\begin{aligned}
&u_i\geq 0, 1\leq i\leq 4,\\
&2u_1+3u_2+6u_3+10u_4\leq 1,\\
&(2+15)u_1+(3+15)u_2+(6-15)u_3+(10-15)u_4\leq 1.\\
\end{aligned}
\right.
\end{equation*}
Separate it into following two systems:
\begin{equation}\label{sysex1}
\left\{
\begin{aligned}
&u_i\geq 0,\ 1\leq i\leq 4,\\
&2u_1+3u_2+6u_3+10u_4\leq 1,\\
&u_1+u_2<u_3+u_4.\\
\end{aligned}
\right.
\end{equation}
and
\begin{equation}\label{sysex2}
\left\{
\begin{aligned}
&u_i\geq 0,\ 1\leq i\leq 4,\\
&(2+15)u_1+(3+15)u_2+(6-15)u_3+(10-15)u_4\leq 1,\\
&u_1+u_2\ge u_3+u_4.\\
\end{aligned}
\right.
\end{equation}
The volume of ($\ref{sysex1}$) is $V_{(2,3),(6,10)}$, and the volume of  ($\ref{sysex2}$) is $V_{(-9,-5),(17,18)}$. Both of them are $(2,2)$-forms.  Moreover, by (\ref{relation2}), we have
\begin{equation*}
    \begin{aligned}
        V_{(2,3),(6,10)}&=V_{(4),(6,8,9)},\\
        V_{(-9,-5),(17,18)}&=V_{(1),(17,8,12)},\\
    \end{aligned}
\end{equation*}
which reduce the calculations of $(2,2)$-forms to the calculations of $(1,3)$-forms.
\end{example*}
\vspace{0.7cm}

On the other hand, by the definition of $(a,b)$-form, we have
\begin{equation}\label{relation1}
    V_{(r_1,\ldots,r_l),(s_1,\ldots,s_{m-1})}+V_{(s_1,\ldots,s_{m-1}),(r_1,\ldots,r_l)}=\frac{1}{(l+m-1)!\prod\limits_{j=1}^l r_j\prod\limits_{i=1}^{m-1} s_i}.
\end{equation}
which means the computation of an $(m-2,l+1)$-form is equivalent to that of an $(l+1,m-2)$-form. 
Combining (\ref{relation2}) and (\ref{relation1}), we have
\begin{equation*}
    V_{(r_1,\ldots,r_l),(s_1,\ldots,s_{m-1})}=\frac{1}{(l+m-1)!\cdot s_1\prod\limits_{i=1}^l(r_i+s_1)\prod\limits_{j=2}^{m-1}(s_j-s_1)}-V_{(s_1,r_1+s_1,\ldots,r_l+s_1),(s_2-s_1,\ldots,s_{m-1}-s_1)}.
\end{equation*}

Note that
\begin{equation*}
    V_{(z_1,\ldots,z_{l+m-1}),\emptyset}=0.
\end{equation*}

Therefore, by induction, 
\begin{equation*}
\begin{aligned}
    V_1=&V_{(r_1,\ldots,r_l),(s_1,\ldots,s_{m-1})}\\
    &=\frac{1}{(l+m-1)!}\sum\limits_{j=1}^{m-1}\frac{(-1)^{j-1}}{s_j\prod\limits_{i=j+1}^{m-1}(s_i-s_j)\prod\limits_{t=1}^{j-1}(s_i-s_t)\prod\limits_{k=1}^l(r_k+s_j)}.
\end{aligned}
\end{equation*}

Similarly,
\begin{equation*}
    V_2=\frac{1}{(l+m-1)!}\sum\limits_{i=1}^l\frac{(-1)^{i-1}}{(r_i+s_m)\prod\limits_{j=i+1}^l(r_j-r_i) \prod\limits_{t=1}^{i-1}(r_i-r_t)\prod\limits_{k=1}^{m-1}(s_k+r_i)}.
\end{equation*}

\vspace{0.7cm}

\begin{example*}
        Let $R=\{2,3,6\},S=\{10,15\}$. Then by (\ref{relation2}) and (\ref{relation1}),
        \begin{equation*}
            \begin{aligned}
                V_{(2,3),(6,10)}&=V_{(4),(6,8,9)}\\
                &=\frac{1}{4!\cdot 4\cdot 6\cdot 8\cdot 9}-V_{(6,8,9),(4)}\\
                &=\frac{1}{4!\cdot 4\cdot 6\cdot 8\cdot 9}-V_{\emptyset,(4,10,12,13)}\\
                &=\frac{1}{4!\cdot 4\cdot 6\cdot 8\cdot 9}-\frac{1}{4!\cdot 4\cdot 10\cdot 12\cdot 13}+V_{(4,10,12,13),\emptyset}\\
                &=\frac{1}{4!\cdot 4\cdot 6\cdot 8\cdot 9}-\frac{1}{4!\cdot 4\cdot 10\cdot 12\cdot 13}.
            \end{aligned}
        \end{equation*}
One can check that this result equals to the formula of $V_1$ given above.
        
    \end{example*}
\vspace{0.7cm}

Now let us compute $V$, the volume of (\ref{sys3}). We have
\begin{equation*}
    \begin{aligned}
        V&=\frac{1}{(l+m-1)!\prod\limits_{j=1}^l r_j\prod\limits_{i=1}^{m-1} s_i}-(V_1+V_2)\\
        &=\left(\frac{1}{(l+m-1)!\prod\limits_{j=1}^l r_j\prod\limits_{i=1}^{m-1} s_i}-V_1\right)-V_2\\
        &=\frac{1}{(l+m-1)!}\sum\limits_{i=1}^l\frac{(-1)^{i-1}}{r_i\prod\limits_{j=i+1}^l(r_j-r_i) \prod\limits_{t=1}^{i-1}(r_i-r_t)\prod\limits_{k=1}^{m-1}(s_k+r_i)}-V_2\\
        &=\frac{1}{(l+m-1)!}\sum\limits_{i=1}^l\frac{(-1)^{i-1}s_m}{r_i\prod\limits_{j=i+1}^l(r_j-r_i)\prod\limits_{t=1}^{i-1}(r_i-r_t)\prod\limits_{k=1}^m(s_k+r_i)}.\\
    \end{aligned}
\end{equation*}
Here the third equality is due to (\ref{relation1}). Therefore,
\begin{equation}\label{eq2}
\begin{aligned}
   p_{R>S}(n)&=\frac{1}{s_m}\cdot V\cdot n^{l+m-1}+o(n^{l+m-1})  \\
   &=A_{l+m-1}\sum\limits_{i=1}^l\frac{(-1)^{i-1}}{r_i\prod\limits_{j=i+1}^l(r_j-r_i)\prod\limits_{t=1}^{i-1}(r_i-r_t)\prod\limits_{k=1}^m(s_k+r_i)}+o(n^{l+m-1}).
\end{aligned}
\end{equation}

Then (\ref{eq3}) comes directly from (\ref{eq1}) and (\ref{eq2}), which completes the proof of Theorem \ref{thm1}.

\section{Discussion of the applications of  Theorem \ref{thm1} to certain $R,S$}\label{sec3}
Let $R_N=\{r,r+m,\ldots,r+m(N-1)\},S_N=\{s,s+m,\ldots,s+m(N-1)\}, I_N=\left[\max(r,s)+m(N-1)\right]\setminus(R_N\cup S_N)$. Here $r,s,m$ are three positive integers such that $r\not\equiv s\ (mod\ m)$ and $(r,s,m)=1$. Denote $\dfrac{p_{R_N>S_N,I_N}(n)}{p_{R_NS_NI_N}(n)}$ by $C_{n,N}$. Then by Theorem \ref{thm1},
\begin{equation*}
    \lim\limits_{n\rightarrow \infty}C_{n,N}=C\cdot\left(\frac{r}{m}+N-1\right)_N\cdot\left(\frac{s}{m}+N-1\right)_N,
\end{equation*}
where 
\begin{equation*}
    C=\sum\limits_{i=1}^N\frac{(-1)^{i-1}}{(N-i)!(i-1)!\
    \left(\frac{s+r}{m}+N+i-2\right)_N\left(\frac{r}{m}+i-1\right)}.
\end{equation*}
Here $(a)_N=\prod\limits_{i=0}^{N-1}(a-i)$.

Let $g(x)=\sum\limits_{i=1}^N\dfrac{(-x)^{i-1}}{(N-i)!(i-1)!\left(\frac{s+r}{m}+N+i-2\right)_N\left(\frac{r}{m}+i-1\right)}$. Note that $C=g(1)$. We have 
\begin{equation}\label{eq4}
    \left[g(x)x^{N+\frac{s+r}{m}-1}\right]^{(N)}=\sum\limits_{i=1}^N\frac{(-1)^{i-1}x^{\frac{s+r}{m}+i-2}}{(N-i)!(i-1)!\left(\frac{r}{m}+i-1\right)}.
\end{equation}

Let $h(x)=\sum\limits_{i=1}^N\dfrac{(-1)^{i-1}x^{\frac{r}{m}+i-1}}{(N-i)!(i-1)!(\frac{r}{m}+i-1)}$. Then 
\begin{equation*}
\begin{aligned}
     h^\prime(x)&=\sum\limits_{i=1}^N\frac{(-1)^{i-1}x^{\frac{r}{m}+i-2}}{(N-i)!(i-1)!}\\
     &=\frac{x^{\frac{r}{m}-1}}{(N-1)!}\sum\limits_{i=1}^N\binom{N-1}{i-1}(-x)^{i-1}\\
     &=\frac{x^{\frac{r}{m}-1}(1-x)^{N-1}}{(N-1)!}.
\end{aligned}
\end{equation*}

Since $h(0)=0$, we have $h(x)=\displaystyle \int_{0}^x\dfrac{t^{\frac{r}{m}-1}(1-t)^{N-1}}{(N-1)!} dt$. Therefore, 
\begin{equation*}
\begin{aligned}
      C&=g(1)\\
      &=g(1)\cdot 1^{N+\frac{s+r}{m}-1}-g(0)\cdot 0^{N+\frac{s+r}{m}-1} \\
      &=\int_{0}^1dx_1\int_{0}^{x_1}dx_2\ldots\int_{0}^{x_{N-1}} dx_{N}\ \left[g(x_N)x_N^{N+\frac{s+r}{m}-1}\right]^{(N)}\\
      &=\int_{0}^1dx_1\int_{0}^{x_1}dx_2\ldots\int_{0}^{x_{N-1}} dx_{N}\ x_N^{\frac{s}{m}-1}h(x_N)\\
      &=\int_0^1 dx_{N}\int_{x_{N}}^1 dx_{N-1}\ldots\int_{x_2}^1 dx_1\ x_N^{\frac{s}{m}-1}h(x_N)\\
      &=\int_0^1 dx_{N}\ x_N^{\frac{s}{m}-1}h(x_N)\cdot \frac{(1-x_N)^{N-1}}{(N-1)!}\\
      &=\frac{1}{(N-1)!(N-1)!}\int_{0}^1\ dx_N\ x_N^{\frac{s}{m}-1}(1-x_N)^{N-1}\int_0^{x_N}\ dt\  t^{\frac{r}{m}-1}(1-t)^{N-1}.\\
\end{aligned}
\end{equation*}

Therefore,
\begin{equation}\label{eq5}
     \lim\limits_{n\rightarrow \infty}C_{n,N}=\dfrac{\left(\frac{r}{m}+N-1\right)_N\cdot\left(\frac{s}{m}+N-1\right)_N}{(N-1)!(N-1)!}\int_{0}^1\ dx\ x^{\frac{s}{m}-1}(1-x)^{N-1}\int_0^x\ dt\  t^{\frac{r}{m}-1}(1-t)^{N-1}.
\end{equation}

\begin{remark*}
By (\ref{eq5}), one can easily show that when $r<s$, $p_{R_N>S_N,I_N}(n)>p_{R_N<S_N,I_N}(n)$ for $n$ large enough .    
\end{remark*}

This integral is actually a generalization of Euler integral. Let $B(p,q)$ denote the Beta function and $\Gamma(x)$ denote the Gamma function. When $s=m$, we have
\begin{equation*}
    \begin{aligned}
        C&=\frac{1}{(N-1)!(N-1)!}\int_{0}^1\ dx\ (1-x)^{N-1}\int_0^x\ dt\  t^{\frac{r}{m}-1}(1-t)^{N-1}\\
        &=\frac{1}{(N-1)!(N-1)!}\int_{0}^1\ dt\  t^{\frac{r}{m}-1}(1-t)^{N-1}\int_t^1\ dx\ (1-x)^{N-1}\\
        &=\frac{1}{N!(N-1)!}\int_{0}^1\ dt\  t^{\frac{r}{m}-1}(1-t)^{2N-1}\\
        &=\frac{1}{N!(N-1)!}B\left(\frac{r}{m},2N\right).
    \end{aligned}
\end{equation*}
Moreover, for $a\geq N-1$,
\begin{equation*}
    (a)_N=\frac{\Gamma(a+1)}{\Gamma(a-N+1)}.
\end{equation*}

Then we could obtain the asymptotic value of $\lim\limits_{n\rightarrow \infty}C_{n,N}$ as $N\rightarrow \infty$ when $s=m$:
\begin{lemma}
Let $s=m$ and $C_{n,N}$ be defined above. Then
\begin{equation*}\label{lem3}
    \lim\limits_{N\rightarrow \infty}\lim\limits_{n\rightarrow \infty}C_{n,N}=\frac{1}{2^{\frac{r}{m}}}.
\end{equation*}
\end{lemma}
\begin{proof}
Recall the relation between Beta function and Gamma function
\begin{equation*}
    B(p,q)=\frac{\Gamma(p)\Gamma(q)}{\Gamma(p+q)},
\end{equation*}
and the Stirling's approximation
\begin{equation*}
    \Gamma(x+1)\sim\sqrt{2\pi x}(\dfrac{x}{e})^x.
\end{equation*} 
We have
\begin{equation*}
\begin{aligned}
     \lim\limits_{n\rightarrow \infty}C_{n,N}&=\frac{\Gamma(\frac{r}{m}+N)\Gamma(\frac{s}{m}+N)B(\frac{r}{m},2N)}{\Gamma(\frac{r}{m})\Gamma(\frac{s}{m})\Gamma(N+1)\Gamma(N)}\\
     &=\dfrac{\Gamma(\frac{r}{m}+N)\Gamma(2N)}{\Gamma(N)\Gamma(\frac{r}{m}+2N)}\\
     &=\dfrac{\Gamma(\frac{r}{m}+N+1)\Gamma(2N+1)}{\Gamma(N+1)\Gamma(\frac{r}{m}+2N+1)}\cdot\frac{N(2N+\frac{r}{m})}{2N(N+\frac{r}{m})}\\
     &\sim\frac{\sqrt{2\pi \left(N+\frac{r}{m}\right)}\left(\dfrac{N+\frac{r}{m}}{e}\right)^{N+\frac{r}{m}}\sqrt{4\pi N}\left(\dfrac{2N}{e}\right)^{2N}}{\sqrt{2\pi N}\left(\dfrac{N}{e}\right)^N\sqrt{2\pi \left(2N+\frac{r}{m}\right)}\left(\dfrac{2N+\frac{r}{m}}{e}\right)^{2N+\frac{r}{m}}}\\
     &\sim\dfrac{\left(N+\frac{r}{m}\right)^{N+\frac{r}{m}}\left(2N\right)^{2N}}{N^N\left(2N+\frac{r}{m}\right)^{2N+\frac{r}{m}}}.
\end{aligned}
\end{equation*}

Note that for fixed $a$, $(x+a)^{x}\sim e^a x^x$ as $x\rightarrow\infty$. Therefore,
\begin{equation*}
 \begin{aligned}
     \frac{(N+\frac{r}{m})^{N+\frac{r}{m}}}{N^N}\cdot\frac{(2N)^{2N}}{(2N+\frac{r}{m})^{2N+\frac{r}{m}}}&\sim\frac{e^\frac{r}{m}(N+\frac{r}{m})^\frac{r}{m}}{e^\frac{r}{m}(2N+\frac{r}{m})^\frac{r}{m}}\\
     &\sim \frac{1}{2^\frac{r}{m}}.
 \end{aligned}   
\end{equation*}
This finishes the proof of lemma \ref{lem3}.
\end{proof}

Now let us look into another double limit of $C_{n,N}$. Let $R=\{r,r+m,\ldots,r+im,\ldots\},S=\{s,s+m,\ldots,s+im,\ldots\},I=\mathbb{Z}^+\setminus(R\cup S)$. Note that for a fixed $n$,  $C_{n,N}=\dfrac{p_{R>S,I}(n)}{p_{RSI}(n)}$ when $N$ is large enough. Therefore, by (\ref{1mm}), for $(r,s,m)=(1,m,m)$,
\begin{equation*}
\begin{aligned}
        \lim\limits_{n\rightarrow \infty}\lim\limits_{N\rightarrow \infty}C_{n,N}&=\frac{1}{2^{\frac{1}{m}}}\\
        &=\lim\limits_{N\rightarrow \infty}\lim\limits_{n\rightarrow \infty}C_{n,N}.
\end{aligned}
\end{equation*}

This fact indicates that the two double limits of $C_{n,N}$ are equal when $n,N$ tend to infinity for $(r,s,m)=(1,m,m)$. Although we could not provide a direct proof of this commutativity, we conjecture that it holds in more general condition:
~\\
\begin{conje}\label{conj1}
Let $r,s,m$ and other notation be defined above. Let  $R=\{r,r+m,\ldots,r+im,\ldots\},S=\{s,s+m,\ldots,s+im,\ldots\}, I=\mathbb{Z}^+\setminus(R\cup S)$. Then
\begin{equation*}
    \begin{aligned}
        \lim\limits_{n\rightarrow \infty}\frac{p_{R>S,I}(n)}{p_{RSI}(n)}
        &=\lim\limits_{n\rightarrow \infty}\lim\limits_{N\rightarrow \infty}C_{n,N}\\
        &=\lim\limits_{N\rightarrow \infty}\lim\limits_{n\rightarrow \infty}C_{n,N}\\
&=\lim\limits_{N\rightarrow\infty}\frac{M\cdot\left(\frac{r}{m}+N-1\right)_N\cdot\left(\frac{s}{m}+N-1\right)_N}{(N-1)!(N-1)!},
\end{aligned}
\end{equation*}
where
\begin{equation*}
    M=\int_{0}^1\ dx\ x^{\frac{s}{m}-1}(1-x)^{N-1}\int_0^{x}\ dt\  t^{\frac{r}{m}-1}(1-t)^{N-1}.
\end{equation*}
Especially, for $s=m$,
\begin{equation*}
    \lim\limits_{n\rightarrow \infty}\lim\limits_{N\rightarrow \infty}C_{n,N}=\lim\limits_{N\rightarrow \infty}\lim\limits_{n\rightarrow \infty}C_{n,N}=\frac{1}{2^{\frac{r}{m}}}.
\end{equation*}
\end{conje}

This conjecture actually describes the commutativity of the following two operations in the computation of the double limits: let $N\rightarrow\infty$ to use finite sets $R_N,S_N,I_N$ to approximate infinite sets $R,S,I$; let $n\rightarrow\infty$ to obtain the asymptotic ratio between two partition numbers.

Here we provide some possible viewpoints that may help to prove this conjecture:

1. Try to make use of the special case of $(r,s,m)=(1,m,m)$ in which the conjecture holds to study the convergence of $C_{n,N}$;

2. Try to apply the discrete versions of some convergence theorems, for example, the Fatou lemma. The fact that $R,S$ are interleaved sequence may help to this method;

3. Try to show that $C_{n,N}$ uniformly converges for $n$ or $N$.

\section{Concluding remarks}\label{sec4}
We finish this paper with some remarks:

1. The necessary condition of Theorem \ref{thm1} is that $R,S,I$ are finite. Otherwise we could not establish the equality between a partition number and the number of integer points inside a certain finite dimensional polytope. Here we conjecture that Theorem \ref{thm1} holds for infinite $I$. That is, the formula (\ref{eq3})  also applies to the cases of unrestricted partitions.

2. As mentioned, the formula in Conjecture \ref{conj1} is a generalization of Euler integral. The study on the asymptotic of (\ref{eq5}) as $N\rightarrow\infty$ when $r,s\neq m$ may be useful to some extent.

3. It would be very interesting to discover some certain cases that the RHS of (\ref{eq3}) is greater than $\dfrac{1}{2}$, which may help to the study of the biases in partitions. For instance, it is natural to conjecture that the RHS of (\ref{eq3}) is greater than $\dfrac{1}{2}$ when $l=m,\ r_i<s_i$ for $1\leq i\leq l$.

\section*{Acknowledgment}
The authors are grateful to the reviewers and editors for their insightful comments, which have significantly improved the manuscript. This work is supported by the National Key Research and Development Program under Grant 2022YFA1004900 and the National Science Foundation of China under Grant 12031011.

\end{document}

\endinput